\documentclass{article}
\usepackage[utf8]{inputenc}
\usepackage{amssymb}
\input xy
\xyoption{all}

\title{Ulrich  and aCM bundles from invariant theory}
\author{Laurent Manivel}

\newtheorem{theorem}{Theorem}[section]
\newtheorem{lemma}[theorem]{Lemma}
\newtheorem{proposition}[theorem]{Proposition}
\newtheorem{corollary}[theorem]{Corollary}

\def\PP{{\mathbf P}}\def\CC{{\mathbf C}}\def\OO{{\mathbf O}}\def\SS{{\mathbf S}}\def\ZZ{{\mathbf Z}}
\def\fg{{\mathfrak g}}\def\fsl{\mathfrak{sl}}\def\fso{\mathfrak{so}}
\def\fgl{\mathfrak{gl}}
\def\cE{\mathcal{E}}\def\cO{\mathcal{O}}\def\cG{\mathcal{G}}\def\cM{\mathcal{M}}
\def\cI{\mathcal{I}}
\def\ra{\rightarrow}\def\lra{\longrightarrow}

\def\qed{\Box}

\begin{document}

\maketitle

\begin{abstract}
We use certain special prehomogeneous representations of algebraic groups in order to 
construct aCM vector bundles, possibly Ulrich, on certain families of hypersurfaces. 
Among other results, we show that a general cubic hypersurface of dimension seven admits
an indecomposable Ulrich bundle of rank nine, and that a general cubic fourfold admits 
an unsplit aCM bundle of rank six. 
\end{abstract}

\section{Introduction} 

Arithmetically Cohen-Macaulay (aCM) vector bundles have been thoroughly investigated since
Horrock's seminal result that on projective spaces, they must be sums of line bundles. On 
quadrics, spinor bundles are non split aCM vector bundles whose rank grows exponentially with
the dimension of the base, and it was conjectured by Buchweitz, Greuel and Schreyer that 
one cannot expect anything better for higher degree hypersurfaces, independently of their
degrees. In particular low rank aCM vector bundles should not exist on hypersurfaces of
large enough dimension, and this has been established for rank two or three. Actually 
only few examples are known, apart from surfaces and threefolds for which specific techniques 
like the Serre construction are available. 

In this note we suggest a construction of aCM bundles on hypersurfaces using techniques 
from invariant theory. The following situation is quite common: an algebraic group $G$ acts 
linearly on a vector space $V$, and this action has an open orbit; moreover the complement of
this open orbit is an irreducible hypersurface $H$. Sometimes, one can construct on $\PP V$ 
a $G$-invariant morphism $A\otimes\cO_{\PP V}(-k)\stackrel{\varphi}{\lra} B\otimes\cO_{\PP V}$, 
for some integer $k$ and some $G$-modules $A$ and $B$. If this morphism is generically an 
isomorphism, then its cokernel $\cE$ will be a $G$-equivariant aCM sheaf supported on the 
hypersurface $H$. 

In order to get an aCM vector bundle, we would need $H$ to be the schematic support of $\cE$, 
which is not always the case even at the generic point of $H$, as we will see. Moreover, the 
hypersurface $H$ is often singular, while we are mainly concerned by smooth hypersurfaces. 
If the singular locus has large codimension, we can avoid it by restricting to suitable 
projective subspaces, and obtain aCM vector bundles on smooth hypersurfaces. Notably, we
will construct aCM vector bundles of rank six on the general quartic fourfold, and interesting
families of aCM bundles of rank three on the general quartic surface. 

Sometimes our construction works with $k=1$ and the upshot is that we get Ulrich bundles 
on certain hypersurfaces, a very special type of aCM bundles that recently attracted 
considerable interest (see \cite{beauville2} for an introduction). We will start this note 
with a construction of Ulrich bundles on cubic hypersurfaces.

\medskip\noindent {\it Acknowledgements}. We would like to thank Daniele Faenzi, Rosa Mir\'o-Roig, 
Joan Pons-Llopis and Marcello Bernardara for useful discussions. 

\section{Ulrich bundles}

\subsection{General construction}

\noindent {\bf Ansatz A}. Suppose given an affine algebraic group $G$ and 
\begin{enumerate}
\item three $G$-modules $V, A, B$, with $A$ and $B$ of the same dimension, 
\item a $G$-equivariant linear map $\varphi : V  \rightarrow Hom(A,B)$ such that $\varphi (v)$ is injective when $v$ is generic.
\end{enumerate}
Then there is an exact sequence on $\PP V$:
$$0\lra A\otimes\cO_{\PP V}(-1)\stackrel{\varphi}{\lra} B\otimes\cO_{\PP V}\lra \cE\lra 0,$$ 
where the reduced support of the sheaf $\cE$ is a $G$-invariant hypersurface $H$. Let us denote 
by $ \iota$ the embedding of $H$ in $\PP V$.
If this hypersurface is irreducible, and if $\cE=\iota_*E$ is the pushforward of a vector bundle 
on $H$, the degree of $H$ and the rank
of $E$ are related by the formula
\begin{equation}\label{fundamental}
\dim A = \dim B = \mathrm{rank} (E)\times \mathrm{deg} (H).
\end{equation}
            Moreover,   the equation $h$ of $H$ verifies, after normalization,
                         $\det (\varphi)=h^r$ if $r=\mathrm{rank} (E)$. 
The vector bundle $E$ on $H$ is then an {\it Ulrich bundle} (see \cite{beauville2} for more general
definitions and other characterizations of Ulrich bundles).

When $r=1$,  $H$ is a {\it determinantal hypersurface}; if $r=2$ and $\varphi$ is skew-symmetric, then $H$ is a {\it Pfaffian hypersurface} (see \cite{beauville}). Good candidates for constructing examples of higher 
ranks are provided by prehomogeneous spaces, even better, by representations with finitely many orbits. 

Note that it may happen, and we will meet examples of this phenomenon, that 
$\cE$ is not the pushforward of a vector bundle, or even a sheaf on $H$. This would 
mean that $\cI_H\cE=0$, which has no reason to happen in general. But there is a canonical filtration 
of $\cE$ by the subsheaves $\cI_H^p\cE$, such that the successive quotients $\cI_H^p\cE/\cI_H^{p+1}\cE
=i_*E_p$ are indeed pushforwards of sheaves on $H$. The fundamental equation (\ref{fundamental}) becomes
\begin{equation}\label{fundamentalbis}
\dim A = \dim B = (\sum_{p\ge 0}(p+1)\mathrm{rank} (E_p))\times \mathrm{deg} (H).
\end{equation}

\noindent {\it Remark.} Observe in particular that if 
\begin{equation}
\dim A = \dim B = \mathrm{rank} (\cE\otimes_{\cO_\PP V}\cO_H)\times \mathrm{deg} (H),
\end{equation}
then $\cE$ must be the pushforward of a vector bundle at least on an open subset of $H$. 
Moreover this open subset is obviously $G$-stable. 

\subsection{Representations with few orbits} 

Suppose that the $G$-module $V$, or its projectivization $\PP V$, has finitely many orbits. 
In particular it has at most one $G$-invariant hypersurface $H$, whose degree $d$ is known. 
If we can find $\varphi : V  \rightarrow Hom(A,B)$ as above, we can immediately deduce the
value of $r=\mathrm{rank} (\cE)$, at least in the case where it is really the pushforward of
a vector bundle $E$ on $H$.

\smallskip
If the hypersurface $H$ is smooth, the vector bundle $E$
is called an {\it Ulrich bundle}. It is known that any smooth hypersurface carries an Ulrich bundle, but the rank may be very high; the proof of \cite{BGS} gives an estimate of the rank which is exponential in the number of monomials in the equation of the 
hypersurface. Moreover this doesn't seem to be an artefact of the proof but an intrinsic difficulty; for example the minimal 
rank of an Ulrich bundle on a smooth hyperquadric is exponential
in the dimension, and it is conjectured  that aCM bundles of smaller ranks do not exist 
on a general hypersurface of any degree $d\ge 2$ \cite{BGS}. 
It remains therefore an interesting question 
to construct Ulrich bundles of rank as small as possible. 

\smallskip
In our previous setting, the hypersurface $H$ is almost never smooth, but 
we can still obtain Ulrich bundles on smooth linear sections of $H$. These smooth 
sections must have dimension small enough so that they do not meet the singular locus of $H$. 
This singular locus of $H$ is of course $G$-invariant, and in the situation where the 
ambient space $\PP V$ has only few orbits, we can hope that its codimension 
in $H$ is not too small.  

\subsection{Examples}

\subsubsection{The Severi series} Suppose that $S\subset\PP V$ is one of the four Severi varieties, of dimension $2a$ with $a=1,2,4,8$. (We refer to \cite{LMfreud} for the geometry of the Severi varieties.) 
The automorphism group of $S$ is reductive and a finite cover $G$ acts linearly on $V$, a $G$-module of dimension $3a+3$. Moreover $G$ has only two non trivial orbit closures, the cone over $S$ and the cone
over its secant variety, which is a cubic hypersurface $H\subset 
\PP V$. In particular the singular locus of $H$ is $S$, and has 
codimension $a+1$. A general linear section of $H$ of dimension 
$d\le a$ will therefore be smooth. 

The list of Severi varieties is the following:
$$\begin{array}{cccccc}
a & & 1 & 2 & 4 & 8 \\
G & & SL_3 & SL_3^2 & SL_6 & E_6 \\
S & & v_2\PP^2 & (\PP^2)^2 & G(2,6) & \OO\PP^2 \\
V & & S^2\CC^3 & (\CC^3)^{\otimes 2} & \wedge^2\CC^6 & V_{27}
\end{array}$$ 

Choose an equation $h$ of $H$, an element of $Sym^3V^*$. By polarization we get a map 
$$\phi : V \longrightarrow Sym^2V^*\subset Hom(V,V^*),$$
equivariant with respect to the semisimple part of $G$. 
Hence on $\PP V$ an exact sequence 
$$0\rightarrow V\otimes \mathcal{O}_{\PP V}(-1)\rightarrow
V^*\otimes \mathcal{O}_{\PP V}\rightarrow\mathcal{E} \rightarrow 0.$$

\begin{lemma}
The sheaf $\mathcal{E}$ is, on the smooth locus of the cubic hypersurface $H$, 
the pushforward of a vector bundle of rank $a+1$.
\end{lemma}

For $a=2$, we deduce that any smooth cubic surface (which is automatically determinantal)
admits an Ulrich bundle of rank three. For $a=4$, we get that any smooth Pfaffian cubic fourfold 
admits an Ulrich bundle of rank five. This is also the case of any smooth cubic threefold, 
since it is automatically Pfaffian. 

For $a=8$, the Severi variety $S$ is the so-called Cayley plane, 
whose automorphism group has type $E_6$. Moreover $V$ is the minimal representation of $E_6$, of dimension $27$, and the invariant cubic in $\PP V$ is called the Cayley hypersurface. 
We deduce:

\begin{proposition}
A smooth linear section of the Cayley hypersurface, of dimension at most eight, is a smooth cubic supporting an Ulrich bundle of rank nine.  
\end{proposition}

Note that cubic eightfolds obtained as linear sections of the Cayley cubic are far 
from being generic; in fact, remarkably they are rational. But linear sections of
dimension seven are general, and we deduce:

\begin{corollary}
A general cubic hypersurface of dimension at most seven supports an 
indecomposable Ulrich bundle of rank nine. 
\end{corollary}

\noindent {\it Remarks}. \begin{enumerate}
\item Ulrich bundles on cubic surfaces and cubic threefolds have been extensively studied in \cite{CH}. 
\item As a consequence of a remarkable numerical coincidence, a general cubic sevenfold $X$ can
be represented as a linear section of the Cayley cubic in {\it finitely many ways} (up to isomorphism).
Each of these representations allows to define an Ulrich bundle $E_X$ on the cubic, that we call a Cayley
bundle. It would be 
interesting to know how many such representations do exist. Is there a unique one?  
\end{enumerate} 

\medskip
The Cayley bundle $E_X$ on the cubic sevenfold $X$ is indecomposable since it is simple by 
\cite[Theorem 3.4]{IMcubics}. It can be described as follows. Recall that the dual variety of the
Severi variety $S^*
\subset \PP V^*$ is the invariant cubic $H\subset \PP V$. In particular there is an 
incidence correspondence
$$\xymatrix{
 & I:=\PP N_{S^*/\PP V^*}  \ar@{->}[rd]   \ar@{->}[ld] \\ \PP V^*\supset S^* & & H\subset \PP V
 }$$
which is an instance of a Kempf collapsing \cite{kempf}. If $X=H\cap \PP L$ 
is a smooth cubic sevenfold, so that $\PP L$ does not meet $S=Sing(H)$, then $X$
is isomorphic to its preimage in $I$, and we therefore get a map 
$p_X: X\rightarrow S^*$. Then 
$$ p_X^*N_{S^*/\PP V^*}=E_X\oplus\cO_X(-1).$$

\subsubsection{Three-forms in seven variables}
Let $V$ be a seven dimensional vector space. Then the action of $GL(V)$ on $\wedge^3V$ 
has only finitely many orbits, which were first classified by Schouten in 1931. 
The generic stabilizer is of type $G_2$. Moreover there is a unique semi-invariant hypersurface 
$H$, with an equation $h$ of degree seven. 

Each three-form  $\omega\in\wedge^3V$ defines a symmetric morphism
$$\varphi(\omega) : \wedge^2V\longrightarrow\wedge^5V.$$

\begin{lemma}
If $\omega$ does not belong to $H$, then $\varphi(\omega)$ is an isomorphism. If  $\omega$ is generic in $H$, then $\varphi(\omega)$ has corank three.
\end{lemma}

\noindent {\it Proof}. 
We choose a basis $e_1,\ldots ,e_7$ of $V$ and denote $e_i\wedge e_j\wedge e_k$ by $e_{ijk}$.
Explicit representatives of 
$V-H$ and of the open orbit in $H$ are, respectively \cite[section 3]{KW}:
$$\begin{array}{rcl}
\omega_0 & = & e_{123}+e_{456}+e_{147}+e_{257}+e_{367}, \\
\omega_1 & = & e_{123}+e_{456}+e_{147}+e_{257}.
\end{array}$$
A straightforward computation shows that $\phi (\omega_0)=0$ is injective, while 
 $\phi (\omega_1)$ has kernel $\langle e_{15}, e_{24}, e_{14}-e_{25}\rangle $. 
 The claim follows. $\qed$

\medskip As a consequence, the exact sequence 
$$0\lra \wedge^2V\otimes \cO_{\PP (\wedge^3V)}(-1)\stackrel{\varphi}{\lra} \wedge^5V\otimes 
\cO_{\PP (\wedge^3V)}\lra\cE\lra 0$$
defines an aCM sheaf $\cE$ which, by the Remark in section 2.1, is the pushforward of 
a vector bundle $E$ on the open orbit in $H$. 
Moreover $\det\varphi=h^3$ up to normalization. 

\begin{proposition}
There exists a $76$-dimensional family of degree seven smooth surfaces $S\subset\PP^3$ admitting an 
indecomposable Ulrich bundle of rank $3$.
\end{proposition}

\noindent {\it Proof}. The heptic hypersurface $H$ is singular in codimension at least three, 
just because the complement of the open orbit has codimension three. We will therefore define
a smooth heptic surface $S=H\cap \PP$ by cutting $H$ with a  
general three-dimensional linear subspace
$\PP\subset \PP (\wedge^3V)$. Moreover we get a rank three vector bundle $E_S$ on $S$,
whose pushforward $\cE_S$ to $\PP$ is involved in an exact sequence 
$$0\lra \wedge^2V\otimes \cO_\PP(-1)\lra \wedge^5V\otimes \cO_\PP\lra\cE_S\lra 0.$$
We claim that $E_S$ is in general indecomposable. Indeed, if it where not, it would admit a 
rank one factor, and since the map $\varphi(\omega)$ is symmetric, the surface $S$ would be symmetric 
determinantal. But an easy dimension count shows that the family of heptic
symmetric determinantal surfaces in $\PP^3$ has dimension at most $63$, while the dimension of 
our family is $\dim G(4,35)-\dim \fsl_7=124-48=76$. $\qed$

\medskip\noindent {\it Remarks}. \begin{enumerate}
\item 
Rank three indecomposable Ulrich bundles are known to exist on 
determinantal heptic surfaces \cite{KleppeMR}. 
We would guess that our surfaces are not determinantal, but we are unable to prove it. 
\item One can check that the singular locus of $H$ has codimension exactly three; more precisely
it coincides with the closure of the $31$ dimensional orbit $\cO$ in $\wedge^3V$. Indeed, this closure
contains all the smaller orbits, so we just need to check that the general point $\omega_2$ of 
$\cO$ is a singular point of $H$. To see this, we can use the natural desingularization of $H$ 
by the total space of the conormal bundle of $G(4,V_7)$ (recall that $H$ is the cone over the 
projective dual of this Grassmannian). The set-theoretical fiber $F$ of a point in $\cO$ is the set 
of $U\in G(4,V_7)$ such that $\omega_2$ belongs to $\wedge^2U\wedge V_7$. We can let 
$$\omega_2=  e_{123}+e_{456}+e_{147}.$$
Then we check that $F\simeq\PP^1\times \PP^1$, which implies  that $\omega_2$ is a singular 
point of $H$. Explicitely, $\omega_2$ belongs to $\wedge^2U\wedge V_7$ if and only if $U=\langle
e_1, e_4, e_{23}, e_{56}\rangle$ for some non zero vectors $e_{23}\in\langle e_2, e_3\rangle$ 
and $e_{56}\in\langle e_5, e_6\rangle$. Thus we cannot, unfortunately, extend our Ulrich bundles
to smooth heptic threefolds.

\item Given a pair $(U\subset V_7, \; \omega\in \wedge^2U\wedge V_7)$, 
the wedge product by $\omega$ sends $\wedge^2U$ to $ \wedge^4U\wedge V_7\simeq  \wedge^4U\otimes (V_7/U)$. For $\omega$ 
generic, this map is surjective, hence has a rank three kernel which is also the kernel of $\phi(\omega)$, as a 
straightforward computation shows. The dual map is generically injective as a map of vector bundles, hence defines an 
exact sequence of sheaves on $H_{reg}$, if $p$ denote the projection to $G(4,V_7)$:
$$0\lra p^*Q^*\lra p^*(\wedge^2U)^*\lra E\lra 0.$$
\end{enumerate}

\section{Arithmetically Cohen-Macaulay  bundles}

\subsection{General construction}

\noindent {\bf Ansatz B}. Suppose given an affine algebraic group $G$ and 
\begin{enumerate}
\item three $G$-modules $V, A, B$, with $A$ and $B$ of the same dimension, 
\item a $G$-equivariant map $\varphi : V  \rightarrow Hom(A,B)$, polynomial of degree $k$,  
such that $\varphi (v)$ is injective when $v$ is generic.
\end{enumerate}
Then there is an exact sequence on $\PP V$:
$$0\lra A\otimes\cO_{\PP V}(-k)\stackrel{\varphi}{\lra} B\otimes\cO_{\PP V}\lra \cE\lra 0,$$ 
where the reduced support of the sheaf $\cE$ is a $G$-invariant hypersurface $H$. 
If this hypersurface is irreducible, and if $\cE$ is the pushforward $\cE=\iota*E$ 
of a vector bundle on $H$, the degree of $H$ and the rank
of the vector bundle $E$ on $H$ are related by the formula
$$k\times \dim A = k\times\dim B = \mathrm{rank} (E)\times \mathrm{deg} (H).$$
            Moreover,  the equation $h$ of $H$ verifies, after normalization,
                         $\det (\varphi)=h^r$ if $r=\mathrm{rank} (\cE)$. 

For example one can suppose that $V$ admits an invariant of degree $k+2$. By polarization we get 
an induced linear map
$$Sym^kV\lra Sym ^2V^*\lra Hom(V, V^*),$$
which can be considered as a degree $k$ map from $V$ to $Hom(V, V^*)$. If $V$ is a representation 
with finitely many orbits, there is a good chance that $\varphi$ drops rank on a hypersurface with an 
open orbit whose complement has high codimension.

Over a smooth linear section $X$ of $H$, we will then obtain an aCM vector bundle $\cE_X$, 
hopefully indecomposable.

\subsection{The Freudenthal series}
Consider a Freudenthal variety $F\subset\PP W$, of dimension $3a+3$ with $a=1,2,4,8$
\cite{LMfreud}. Again the automorphism group 
of $F$ is reductive and a finite cover $G$ acts linearly on $W$, a $G$-module of dimension $6a+8$. Here $G$ has only 
three non trivial orbit closures, the cone over $F$, the cone
over its tangent variety, which is a quartic hypersurface $H\subset 
\PP W$, and the singular locus $\Omega$ of this hypersurface, whose dimension is $5a+4$. 
In particular, a general linear section of $H$ of dimension $d\le a+1$ will be smooth. 
The series can be extended to $a=0$ and is modeled on the rational normal cubic:

$$\begin{array}{ccccccc}
 a & & 0 & 1 & 2 & 4 & 8 \\
 G & & SL_3^3 & Sp_6 & SL_6 & Spin_{12} & E_7 \\
 F & & (\PP^1)^3 & LG(3,6) & G(3,6) & OG(6,12) & F_{27} \\
 W & & (\CC^2)^{\otimes 3} & \wedge^{\langle 3\rangle}\CC^6 & \wedge^3\CC^6 & 
 \Delta & W_{56}
 \end{array}$$
 
\smallskip

\begin{lemma}
Let $\fg$ denote the Lie algebra of $G$. 
There exists a non zero quadratic map $\theta : Sym^2W\lra \fg$. 
\end{lemma}

\noindent {\it Proof}. 
The action of $G$ on $W$ yields a map $\fg\rightarrow End(W)$. 
Since $W$ has an invariant symplectic form, this map factorizes through $Sym^2W$. 
Since $W$ and $\fg$ (by the Killing form) are both self-dual, we get the desired map. 
$\qed$

\medskip
Note that $h=K\circ \varphi$, where $K$ denotes the Killing form,
is an equation of the invariant quartic. 

\begin{corollary}
For any $\fg$-module $M$ there exists an equivariant quadratic map $\theta_M : W\lra End(M)$.
If $\theta_M(w)$ is invertible for a generic $w$, one therefore gets an aCM sheaf $\cE_M$
whose reduced support is the hypersurface $H$, defined by an exact sequence
$$0\lra M\otimes\cO_{\PP W}(-2)\stackrel{\theta_M}{\lra} M\otimes\cO_{\PP W}\lra \cE_M\lra 0.$$ 
\end{corollary}

A uniform example, in the sense that it exists for the whole Freudenthal series, is
given by $W$ itself. 

\begin{lemma}
Let $w$ be a general element of $W$. Then $\theta_W(w)$ is invertible. 
\end{lemma}

\noindent {\it Proof}. 
We know that $W$ is a fundamental representation of $G$. Let $v_+$ and $v_-$ be a highest and a 
lowest weight vector. They both belong to the cone over $F$, but the sum $v_++v_-$ belongs to the cone over 
the open orbit,  and is therefore a general element of $W$. Since $\fg\subset Sym^2W^*$ is the space 
of quadratic equations of $F$, we have $a(v_++v_-)=2a(v_+,v_-)$. This must be a multiple of the center of the 
reductive group defined as the common stabilizer of the highest and lowest weight vectors, hence a multiple 
of the fundamental coweight $H$ associated to $W$. So the claim reduces to the fact 
that for any weight $\omega$ of $W$, $\omega (h)$ is non zero, which we check case by case.
$\qed$ 

\medskip
Unfortunately, we will see later that the sheaf $\cE_W$ is not a vector bundle on the 
smooth locus of $H$. But we can get genuine aCM bundles from the most basic representations
of $\fsl_6$ and $\fso_{12}$, that will be studied in the next two sections.

\subsection{Quartic surfaces}

Consider the case where $G=SL_6=SL(V_6)$, $F=G(3,V_6)$ and $W=\wedge^3V_6$, where $V_6$ denotes
a six-dimensional complex vector space. 
The map $Sym^2W\ra\fsl_6$ can be described as follows (note that by the Schur lemma, there is a unique such 
equivariant projection, up to constant). There is a natural map from $W\otimes W$ to $\wedge^5V_6\otimes V_6$ 
defined by the composition 
$$\wedge^3V_6\otimes \wedge^3V_6\lra \wedge^3V_6\otimes \wedge^2V_6\otimes V_6\lra \wedge^5V_6\otimes V_6,$$
where the first map is induced by the polarization morphism $\wedge^3V_6\ra \wedge^2V_6\otimes V_6$ and the 
second one by the exterior product $ \wedge^3V_6\otimes \wedge^2V_6\ra  \wedge^5V_6$. 
As an $SL_6$-module, 
$\wedge^5V_6$ is isomorphic to $V_6^*$. So the symmetric part of the previous morphism gives a map
$$Sym^2( \wedge^3V_6)\lra V_6^*\otimes V_6\lra \fsl_6,$$
where the rightmost arrow is the projection to the traceless part. If the resulting morphism is non zero,
it must be the map we are looking for. 

Let us denote by $\theta : \wedge^3V_6\ra \fsl_6$ the corresponding quadratic map, and let us compute it explicitely. 
A general element in $\PP (\wedge^3V_6)$ belongs to a unique secant line to $G(3,V_6)$, so a general element of 
$\wedge^3V_6$ can be written, up to scalar, as $\omega_0=e_{123}+e_{456}$, where $e_{123456}$ is a fixed generator of $\wedge^6V_6$. 
We get 
$$\theta (\omega_0)=Id_{\langle e_1, e_2, e_3\rangle}-Id_{\langle e_4, e_5, e_6\rangle},$$
a semisimple endomorphism whose eigenspaces are precisely the three-spaces $A=\langle e_1, e_2, e_3\rangle$ and 
$B=\langle e_4, e_5, e_6\rangle$. The eigenvalues of $\theta (\omega_0)$ on $\wedge^3V_6$ are $3$ on $\wedge^3A$, 
$1$ on $\wedge^2A\otimes B$, $-1$ on $A\otimes \wedge^2B$ and $-3$ on $\wedge^3B$; in particular $\theta (\omega_0)$
acts on $\wedge^3V_6$ as an isomorphism, confirming the previous lemma. More generally, $\theta_M (\omega_0)$ will be an 
isomorphism for any irreducible module $M=S_\mu V_6$ such that $\mu$ is a partition of {\it odd} size. 

A general element $\omega_1$ of the quartic hypersurface $H\subset\wedge^3V_6$ is a general element of the affine tangent space to 
$G(3,V_6)$ at a uniquely determined point $A=\langle e_1, e_2, e_3\rangle$. Since the affine tangent space to the Grassmannian
at that point is $\wedge^2A \wedge V_6$, we can choose $\omega_1= e_{126}+e_{234}+e_{315}$. Then we get 
$$\theta (\omega_1)=e_4^*\otimes e_1+e_5^*\otimes e_2+e_6^*\otimes e_3,$$
a square zero endomorphism whose kernel and image are precisely the three-space $A$.

If we apply our construction to the natural representation $V_6$ of $G$, we obtain an aCM sheaf 
$\cE_{V_6}$ with an exact sequence 
$$0\rightarrow V_6\otimes \mathcal{O}_{\PP W}(-2)\rightarrow
V_6\otimes \mathcal{O}_{\PP W}\rightarrow\mathcal{E}_{V_6}\rightarrow 0.$$
Our computation of $\theta (\omega_1)$ shows that $\mathcal{E}_{V_6}$ has  rank $3$ on 
the open orbit $\cO$ of $H$ (which is also it smooth locus). Since $2\times 6=3\times 4$, we can
conclude that $\mathcal{E}_{V_6}$ is the pushforward of a vector bundle $E_{V_6}$ on $\cO$. 
An alternative description of this 
bundle is the following: since $H$ is the tangent variety to $F=G(3,V_6)$ (and also its projective  
dual), the open orbit $\cO$ is naturally fibered over $G(3,V_6)$ (with fibers the projectived
tangent spaces; the natural map to $H$ is known to be birational, hence a resolution of 
singularities \cite{LMfreud}). Let us denote by $\pi$ the projection from $H_{reg}$ to $G(3,V_6)$.
From the explicit description of $\theta (\omega_1)$ above, we see that 
$$E_{V_6|H_{reg}}=\pi^*Q, \qquad E_{V_6^{\vee}|H_{reg}}=\pi^*U^*,$$
if we denote by $U$ and $Q$ the tautological and the quotient rank three vector bundles on 
$G(3,V_6)$. We deduce:

\begin{proposition}
Let $S$ be a general quartic surface in $\PP^3$, represented as a linear section of the invariant 
quartic $H\subset\PP(\wedge^3V_6)$. The vector bundles $E_{V_6|S}$ and  $E_{V_6^{\vee}|S}$ are 
indecomposable rank three aCM vector bundles on $S$.
\end{proposition}

\noindent {\it Proof.} We need to check that $E :=E_{V_6|S}$ is indecomposable. 
Suppose it is not. Then it must split as the direct sum $L\oplus M$ of a line bundle $L$ 
and an aCM rank two vector bundle $M$ on $S$. Both are generated by global sections, 
and $h^0(L)+h^0(M)=6$. But in general $Pic(S)=\ZZ$, and $L$ cannot be trivial because 
$h^0(E)=6=\chi(E)$ and $h^1(E)=0$, hence $h^0(E^*)=h^2(E)=0$. So $L=\cO_S(k)$ for some
$k>0$, which implies that $h^0(L)\ge 4$. Therefore $h^0(M)\le 2$, but since $M$ is generated
of rank two this is only possible if $M$ is trivial, which again contradicts $h^0(E^*)=0$.
$\qed$

\medskip We will call the aCM bundles $E$ that we obtain on a quartic surface $S$ by this construction, 
the {\it Grassmann vector bundles} on $S$. The pair $(S,E)$ is defined by a linear 
section of the invariant quartic $H$, so the number of parameters is $\dim G(4,20)-\dim SL_6=4\times 16-35=29$. Since there are $19$ parameters for $S$, we expect that a general quartic surface $S$ 
supports a ten-dimensional family of Grassmann bundles (or two such families).  

\smallskip
Note that the rational projection from $H$ to $G$ is defined by the derivatives
of the quartic invariant, which implies that $c_1(E)=3h$ if $h$ denotes the hyperplane class. 
Moreover the Hilbert polynomial $P_E$ of $E$ is easily computed, and by Riemann-Roch, we get 
$$P_E(k)= 6(k+1)^2 = 6+ch_2(\cE(k)) = 6+ch_2(\cE)+k c_1(\cE)h+k^2\frac{h^2}{2},$$
hence $ch_2(\cE)=0$. 

\begin{proposition}
The two Grassmann bundles $E=\cE_{V_6|S}$ and $F=\cE_{V_6^{\vee}|S}$ on $S$ are in natural
duality, in the sense that $E^*\simeq F(-2)$.
\end{proposition}

\noindent {\it Proof.} If we consider $S$ as a subvariety of $\PP\hat{T}_G$, the line bundle 
$\cO_S(-1)$ is identified with $\cO_{\hat{T}_G}(-1)$ and is therefore a subbundle of $\pi^*\hat{T}_G$. 
Since $T_G(-1)$ is a quotient of $\hat{T}_G$, we get a map from $\cO_S(-1)$ to $\pi^*(T_G(-1))$, which 
is never zero since $S$ does not meet $G$. But $T_G=U^*\otimes Q$, hence a natural injective map
$$\cO_S(-1)\lra E\otimes F(-3).$$
The image of this map is everywhere non degenerate and yields the desired duality. 
$\qed$

\medskip
The Mukai vector of $E$ is $v(E)=(3,3h,3)$, in particular it is not primitive.
Moreover the dimension of the moduli space $\cM$ of simple vector bundles on $S$ at the point defined by $E$ is $v(E)^2+2=9\times 4-2\times 9+2=20$. A simple diagram chasing shows that the kernel of the 
moduli space at the point $[E]$ to the family of Grassmann bundles 
should be given by the exact sequence 
$$0\lra T_{[E]}\cG \lra Hom(L,\wedge^3V_6)/\fgl_6 \lra Sym^4L^*/\langle H_L\rangle \lra 0.$$
Here we denoted by $L$ the four-dimensional subspaces of $\wedge^3V_6$ such that $S=H\cap \PP L$,
and by $H_L$ the restriction to $L$ of an equation of $H$. There is a natural map that sends 
$u\in Hom(L,\wedge^3V_6)$ to the quartic polynomial $H(u(x),x,x,x)$ on $L$; if $u$ comes from 
$\fsl_6$, this polynomial is identically zero by invariance.

\medskip\noindent {\it Fantasy.} We thus get a subvariety of dimension ten in a 
hyperkaehler moduli space of dimension twenty. Is it Lagrangian? 

The tangent space to the moduli space at a point defined by a Grassmann vector bundle is a 
twenty-dimensional vector space endowed with a non degenerate skew-symmetric form. Is there
a natural identification of this tangent space with $(\wedge^3V_6, \wedge)$? 

\medskip\noindent {\it Remark.} 
The action of $\theta (\omega_1)$ on $\wedge^3V_6$ is nilpotent of order four. In terms of $B=\langle e_4, e_5, e_6\rangle$
(which is not uniquely defined by $\omega_1$), this action is graded, in the sense that it maps $\wedge^3B$ to $A\otimes \wedge^2B$
(injectively), $A\otimes \wedge^2B$ to $\wedge^2A\otimes B$ (isomorphically), $\wedge^2A\otimes B$ to $\wedge^3A$ (surjectively),
and $\wedge^3A$ to zero. This implies that the kernel of $\theta (\omega_1)$ is the direct sum of  $\wedge^3A$ with a 
hyperplane in $\wedge^2A\otimes B$, so the cokernel has dimension $11$ and not $20\times 2/4=10$.
This implies that the $\cE_{\wedge^3V_6}$ cannot be the pushforward of a vector bundle on the 
smooth locus of $H$; instead, it is the extension of the pushforward of a vector bundle of rank $10$
by the pushforward of a line bundle, neither of which has is aCM. Nevertheless, $\cE_{\wedge^3V_6}$
is an interesting example of an aCM sheaf supported on the first infinitesimal neighbourhood of $H$. 
See \cite{BHMP} for more on aCM sheaves on non reduced schemes.

\subsection{Quartic fourfolds}

In this section we consider the case where $G=Spin_{12}=Spin(V_{12})$, where $V_{12}$
is a twelve dimensional vector space endowed with a non degenerate quadratic form. 
Here $F=\SS_{12}$ parametrizes one of the two families of maximal isotropic subspaces
of $V_{12}$ and $W=\Delta_+$ is one of the half-spin representations. Recall that 
these representations can be described explicitely by choosing a decomposition of 
$V_{12}$ as the direct sum of two isotropic subspaces $E$ and $E'$. Since the quadratic
form defines a perfect duality between them, we rather write $V_{12}=E\oplus E^\vee$. 

Then the Lie algebra of $G$ decomposes as 
$$\fg=\fso(V_{12})\simeq \wedge^2V_{12}=\wedge^2E\oplus End(E)\oplus \wedge^2E^*.$$
Moreover the half-spin representation maybe defined as
$$W=\Delta_+ = \CC\oplus \wedge^2E\oplus \wedge^4E\oplus \wedge^6E,$$
with the natural action of $\fg$ defined by wedge products and contractions. 

This allows to describe the natural quadratic map $\theta: W\ra\fg$ as follows. Let us 
fix an orientation of $E$, that is, a generator of $\wedge^6E$. It will be convenient to
fix a basis $e_1,\ldots ,e_6$ of $E$ and choose $e_{123456}\in \wedge^6E$ for generator. 
With respect
to the previous decompositions, $\theta$ maps a spinor $w=(\omega_0, \omega_2, \omega_4,\omega_6)$
to $\theta(w)=(\theta_+(w),\theta_0(w),\theta_-(w))$ with 
$$
\theta_+(w) = \omega_6\omega_2-\omega_4*\omega_4, \quad
\theta_0(w)  = \omega_0\omega_6 Id_E-\omega_2*\omega_4 \quad
\theta_-(w)  =  \omega_0\omega_4-\omega_2\wedge\omega_2. 
$$
One will readily check that $\theta$ gives the quadratic equations of the spinor variety
$\SS_{12}$, which is the closure of the image of $\wedge^2E$ to $\PP\Delta_+$ given by
$$\omega \mapsto [1, \omega ,\omega^2 ,\omega^3 ].$$ 

\smallskip Let us compute explicitely $\theta$ on the open orbit $\Delta_+-H$, and on the 
open orbit in $H$. 

Recall that $H$ is the tangent hypersurface to the spinor variety $\SS_{12}$, so 
a general point of $H$ is a general point of the projective tangent space of $\SS_{12}$ at a
general point $p$. Since $\SS_{12}$ is homogeneous, we can choose $p$ to be any point, 
say $p=[1,0,0,0]$. 
It follows from the explicit parametrization of (an open subset of) $\SS_{12}$ given above that 
the corresponding tangent space is given by the set of points of the form 
$w_1=[\omega_0, \omega_2,0,0]$. 
We get a general point if $\omega_0$ is non zero and $\omega_2$ has maximal rank, so we can let 
$\omega_0=1$ and $\omega_2=e_{12}+e_{34}+e_{56}$. We then get, up to a non zero scalar, 
$$\theta(w_1)=(0,0,e^*_{12}+e^*_{34}+e^*_{56})\in\fso (V_{12}).$$
In particular $\theta_{V_{12}}(w_1)$ is a rank six endomorphism whose cokernel is $V_{12}/E\simeq E^*$.  

Now let $w_0=[1,0,0,e_{123456}]$. Then $\theta(w_0)=Id_E$, which, when considered as an endomorphism
of $V_{12}$ yields $$\theta_{V_{12}}(w_0)= Id_E-Id_{E^*}.$$ 
In particular $\theta_{V_{12}}(w_0)$ is invertible (which implies a posteriori that $w_0$ does not
belong to $H$). We thus get a sheaf $\cE_{V_{12}}$ supported on $H$, and an exact sequence 
$$0\rightarrow V_{12}\otimes \mathcal{O}_{\PP \Delta_+}(-2)\rightarrow
V_{12}\otimes \mathcal{O}_{\PP \Delta_+}\rightarrow\mathcal{E}_{V_{12}}\rightarrow 0.$$
Since the rank of the cokernel of $ \theta_{V_{12}}(w_1)$ is $6=12\times 2/4$, we can assert 
that on the smooth part of $H$, $\mathcal{E}_{V_{12}}$ is the pushforward of a vector bundle
$E_{V_{12}}$. 
Moreover, as in the previous case $H_{reg}$ is fibered over $F=\SS_{12}$, and if we denote
the projection by $\pi$ we get that
$$E_{V_{12}}=\pi^*U^*,$$
if $U$ denotes the tautological rank six vector bundle on the spinor variety $\SS_{12}\subset G(6,V_{12})$. Note that as in the Grasmmannian case, the natural map $$\cO(-1)\ra\pi^*\hat{T}
\ra\pi^* T(-1)=\pi^*(\wedge^2U^*(-1))$$ implies that $E_{V_{12}}(-1)$ is self-dual.
We deduce:

\begin{theorem}
The general quartic fourfold admits an unsplit aCM vector bundle of rank six.
\end{theorem}

\noindent {\it Proof.}    
It was checked in \cite[Theorem 2.2.1]{abuaf} that the general quartic
fourfold $X$ is a linear section of the $Spin_{12}$-invariant quartic $H$, 
in finitely many ways. (Note that the number of such representations is not known.)
In particular the 
intersection of $H$ with a general $\PP^5$ does not hit its singular locus, and is therefore
contained in the open orbit. Restricting $E_{V_{12}}$, 
we therefore get a vector bundle $E_X$ on $X$ and an exact sequence
$$0\rightarrow V_{12}\otimes \mathcal{O}_{\PP^5}(-2)\rightarrow
V_{12}\otimes \mathcal{O}_{\PP^5}\rightarrow \iota_* E_X\rightarrow 0.$$
This bundle $E_X$ is unsplit:
if it was a direct sum of line bundles, by self-duality and since $E_X$ is
generated by global sections we would have $E_X\simeq a\cO_X\oplus a\cO_X(2)\oplus (6-2a)\cO_X(1)$,
and then $h^0(E_X)=36+10a\ne 12$, a contradiction.$\qed$ 

\medskip
Let us call those aCM vector bundles on $X$ its {\it spinor bundles}. We figure out that these 
bundles should be indecomposable but we have not been able to prove it. We get one such 
bundle for each of the finitely many representations of $X$ as a linear section of the 
$Spin_{12}$-invariant quartic $H$, and 
we would guess that those bundles cannot be isomorphic one to the others. 
In fact, there is probably a way to reconstruct the embedding of $X$ in the projectivized 
half-spin representation from the spinor bundle, but we do not know how to do that.

\medskip
Institut de Math\'ematiques de Toulouse, UMR 5219, Universit\'e de Toulouse, CNRS, UPS IMT F-31062 Toulouse Cedex 9, France.

\begin{thebibliography}{Aa}


\bibitem{abuaf} Abuaf R., {\it 
    On quartic double fivefolds and the matrix factorizations of exceptional 
    quaternionic representations}, arXiv:1709.05217.
    
\bibitem{AHMP}
Aprodu M., Huh S., Malaspina F., Pons-Llopis J.,  {\it   
    Ulrich bundles on smooth projective varieties of minimal degree}, 
    arXiv:1705.07790.

\bibitem{BHMP}
Ballico E., Huh S., Malaspina F., Pons-Llopis J.,  {\it    
    ACM sheaves on the double plane}, arXiv:1604.00866.

\bibitem{beauville} Beauville A., {\it Determinantal hypersurfaces},     Michigan Math. J.
{\bf 48}  (2000), 39--64.

\bibitem{beauville2} Beauville A., {\it An introduction to Ulrich bundles}, 
arXiv:1610.02771.

\bibitem{BGS} Buchweitz R.-O., Greuel G.-M., Schreyer, F.-O.,
{\it Cohen-Macaulay modules on hypersurface singularities} II, 
Invent. Math. {\bf 88} (1987),  165--182.

\bibitem{CH} Casanellas M., Hartshorne R., {\it Stable Ulrich bundles},
International J. Mathematics {\bf 23} (2012), 1250083--1250133.

\bibitem{coskun}
Coskun I., {\it Ulrich bundles on quartic surfaces with Picard number 1},
C. R. Acad. Sci. Paris, Ser. I {\bf 351} (2013) 221--224.

\bibitem{FPL}
Faenzi D., Pons-Llopis J., 
{\it The Cohen-Macaulay representation type of arithmetically Cohen-Macaulay varieties},
arXiv:1504.03819.

\bibitem{IMcubics} 
Iliev A., Manivel L., {\it On cubic hypersurfaces of dimensions $7$ and $8$}, 
Proc. London Math. Soc. {\bf 108} (2014), 517--540.
	
\bibitem{kempf}
Kempf G., {\it On the collapsing of homogeneous bundles},
Invent. Math. {\bf 37} (1976), 229--239.
	
\bibitem{KleppeMR} Kleppe J., Mir\'o-Roig R., {\it  The representation type of
determinantal varieties},
Algebras and Representation Theory {\bf 20} (2017), 1029--1059.

\bibitem{MRPL} Mir\'o-Roig R., Pons-Llopis J., {\it  N-dimensional Fano varieties of 
wild representation type}, J. Pure Appl Algebra. {\bf 218} (2014), 1867--1884.

\bibitem{KW} Kraskiewicz W., Weyman J., {\it 
Geometry of orbit closures for the representations associated to gradings of Lie algebras of types} $E_7$, arXiv:1301.0720.

\bibitem{LMfreud} Landsberg J., Manivel L.,
{\it The projective geometry of Freudenthal's magic square}, Journal of Algebra {\bf 239} (2001),
477--512.   



\end{thebibliography}
\end{document}